\pdfoutput=1
\documentclass[11pt]{article}
\usepackage[left=1in,right=1in,top=1in,bottom=1in]{geometry}
\usepackage{times}
\usepackage{expl3}
\usepackage{cite}
\usepackage[table]{xcolor}
\usepackage{multirow}
\usepackage{stackengine} 
\usepackage{hhline}
\usepackage{lipsum}
\usepackage{titlesec}
\usepackage{wrapfig}
\usepackage{enumerate}
\usepackage{epsfig}
\usepackage{amsmath}
\usepackage{tabularx}
\usepackage{array}
\usepackage{booktabs}
\usepackage{enumitem}
\usepackage{bbm}
\usepackage{calc}
\usepackage{graphicx}
\usepackage{amsmath}
\usepackage[title]{appendix}
\usepackage{amssymb}
\usepackage{epstopdf}
\usepackage{boldline}
\usepackage{arydshln}
\usepackage{calligra}
\usepackage{bm}
\usepackage{url}
\usepackage{blindtext}
\usepackage{accents}

\newcommand{\define}{\stackrel{\mbox{\tiny def}}{=}}

\newtheorem{definition}{Definition}
\newtheorem{theorem}{Theorem}
\newtheorem{corollary}{Corollary}
\newtheorem{lemma}{Lemma}

\newtheorem{example}{Example}
\newtheorem{remark}{Remark}
\usepackage{mathtools}
\usepackage{epstopdf}
\usepackage{balance}
\usepackage{thmtools}
\usepackage{thm-restate}
\usepackage{hyperref}
\usepackage{cleveref}
\usepackage[mathscr]{euscript}

\usepackage[ruled,vlined]{algorithm2e}
\include{pythonlisting}
\newcommand{\ostar}{\mathbin{\mathpalette\make@circled\star}}

\makeatletter
\newcommand{\removelatexerror}{\let\@latex@error\@gobble}
\makeatother
\makeatletter
\newcommand*{\rom}[1]{\expandafter\@slowromancap\romannumeral #1@}
\makeatother

\ExplSyntaxOn
\newcommand\latinabbrev[1]{
  \peek_meaning:NTF . {
    #1\@}%
  { \peek_catcode:NTF a {
      #1.\@ }%
    {#1.\@}}}
\ExplSyntaxOff

\titleclass{\subsubsubsection}{straight}[\subsubsection]

\begin{document}
\vspace{1cm}
\title{Generalized Choi-Davis-Jensen's Operator Inequalities and Their Applications}\vspace{1.8cm}

\author{Shih~Yu~Chang and Yimin Wei
\thanks{Shih Yu Chang is with the Department of Applied Data Science,
San Jose State University, San Jose, CA, U. S. A. (e-mail: {\tt
shihyu.chang@sjsu.edu}). Yimiin Wei (Corresponding author) is with the School of Mathematical
Sciences and Key Laboratory of Mathematics for Nonlinear Sciences, Fudan University, Shanghai, 200433, P. R. China.
This author is supported by the Shanghai Municipal Science and Technology Commission under grant 23WZ2501400
and the Ministry of Science and Technology of China under grant G2023132005L  (e-mail: {\tt
ymwei@fudan.edu.cn, yimin.wei@gmail.com}). 
           }}

\maketitle

\begin{abstract}
The original Choi-Davis-Jensen's inequality, with its wide-ranging applications in diverse scientific and engineering fields, has motivated researchers to explore generalizations. In this study, we extend Davis-Choi-Jensen’s inequality by considering a nonlinear map instead of a normalized linear map and generalize operator convex function to any continuous function defined in a compact region. The Stone–Weierstrass theorem and Kantorovich function are instrumental in formulating and proving generalized Choi-Davis-Jensen's inequalities. Additionally, we present an application of this generalized inequality in the context of statistical physics.
\end{abstract}


\section{Introduction}\label{sec: Introduction}

Choi~\cite{choi1974schwarz} and Davis~\cite{davis1957schwarz} demonstrated that for a normalized positive linear map $\Phi$ and an operator convex function $f$ defined on an interval $I$, the following inequality known as Choi-Davis-Jensen’s inequality holds for any self-adjoint operator $\bm{A}$: 
\begin{eqnarray}\label{eq: CDJ Inquality}
f(\Phi(\bm{A}))&\leq&\Phi(f(\bm{A})),
\end{eqnarray}
where the operator inequality $\leq$ is defined as Loewner ordering, i.e., $\bm{B} \leq \bm{A}$ if $\bm{A} - \bm{B}$ is a positive-definite operator. Choi-Davis-Jensen’s inequality is the extension of classical Jensen's inequality to the context of operators in functional analysis. Choi-Davis-Jensen's inequality finds applications in various areas, particularly in functional analysis, information theory, quantum mechanics and random operator theory. In functional analysis, Choi-Davis-Jensen's inequality is used to establish norm inequalities for operators acting on Banach spaces by providing bounds on the norm of certain functions of operators~\cite{becker1989functional}. In quantum mechanics, operator Jensen inequality is applied to derive bounds on various entanglement measures. These measures quantify the amount of entanglement in quantum systems~\cite{briet2009properties}. In matrix analysis, Choi-Davis-Jensen's inequality is used to obtain bounds on the spectral radius of matrices. This is essential in the study of eigenvalues and spectral properties of matrices~\cite{briat2011convergence}. In information theory, Choi-Davis-Jensen's inequality plays a role in deriving bounds on the quantum Fisher information, a measure of statistical distinguishability in quantum parameter estimation~\cite{majtey2005jensen}. In statistical machine learning, especially in multivariate analysis, Choi-Davis-Jensen's inequality is used to establish bounds on expectations of certain matrix-valued functions, providing insights into statistical properties~\cite{futami2021loss}. In noncommutative probability, Choi-Davis-Jensen's inequality is employed in the study of noncommutative probability spaces, providing inequalities for the expectation values of noncommutative random variables and the tail bounds of random operators, e.g., matrices and tensors, ensemble~\cite{chang2023algebraic,chang2023tail,chang2023random,chang2023tailRP,chang2022general,chang2022random,chang2022convenient}. In summary, Choi-Davis-Jensen's inequality is a powerful tool in analyzing the properties of operators in various mathematical and physical contexts. They allow for the extension of classical results to the realm of functional analysis, providing insights into the behavior of operators and their applications in diverse fields.

As Choi-Davis-Jensen's inequality has many applications in science and engineering, many authors have tried to generalize this inequality. In~\cite{micic2017choi}, the auhtos present an inequality of Choi-Davis-Jensen type without relying on convexity. Utilizing their principal findings, they establish new inequalities that enhance existing results and demonstrate their generalized Choi-Davis-Jensen inequalities applicable to relative operator entropies and quantum mechanical entropies. Recently, the authors in~\cite{hashemi2021reversing} first provide a better estimate of the second inequality in Hermite-Hadamard inequality and apply this to obtain the reverse of the celebrated Davis-Choi-Jensen’s inequality. Other works in studing the generalization of Choi-Davis-Jensen’s inequality can also be found at references in~\cite{micic2017choi,hashemi2021reversing}. We generalize Davis-Choi-Jensen’s inequality by allowing the function $f$ in Eq.~\eqref{eq: CDJ Inquality} to be continous at a compact region and the mapping $\Phi$ to be a nonlinear map, instead of normalized linear map. Stone–Weierstrass theorem and Kantorovich function are two main ingredits used here to establish generalized Choi-Davis-Jensen's inequality. 

The remainder of this paper is organized as follows. In Section~\ref{sec: Generalized CHOI-DAVIS-JENSEN’S INEQUALITY with error control}, generalized Choi-Davis-Jensen's inequality is derived based on Stone-Weierstrass Theorem and Kantorovich function. The application of the new generalized Choi-Davis-Jensen's inequality to statistical physics is presented in Section~\ref{sec: Applications}. 

\textbf{Nomenclature:}

Inequalities $\geq, >, \leq$  and $<$ besides operators are based on Loewner ordering. The symbol $\Lambda(\bm(A))$ represents the spectrum of the operator $\bm{A}$, i.e., the set of eigenvalues of the operator $\bm{A}$. If $\Lambda(\bm(A))$ are composed by real numbers, we use $\min(\Lambda(\bm(A)))$ and $\max(\Lambda(\bm(A)))$ to represent the minimum and the maximum numbers within $\Lambda(\bm(A))$. Given $M>m>0$ with any $r \in \mathbb{R}$ and $r \neq 1$, Kantorovich function with respect to $m,M,r$ is defined by
\begin{eqnarray}\label{eq: Kantorovich function}
\mathscr{K}(m,M,r)&=&\frac{(mM^r - Mm^r)}{(r-1)(M-m)}\left[\frac{(r-1)(M^r-m^r)}{r(mM^r - Mm^r)}\right]^r.
\end{eqnarray}

\section{Generalized Choi-Davis-Jensen's Inequalities with Error Control}\label{sec: Generalized CHOI-DAVIS-JENSEN’S INEQUALITY with error control}

In this section, we will establish generalized Choi-Davis-Jensen's inequality. Let us recall Stone-Weierstrass Theorem. Consider a continuous real-valued function, denoted as $f(x)$, defined on the closed real interval $[m, M]$, where $m,M \in \mathbb{R}$ and $M>m$. For any given positive $\epsilon$, there exists a polynomial denoted as $p(x)$, such that for all $x$ within the interval $[m, M]$, the absolute difference between $f(x)$ and $p(x)$ is less than $\epsilon$. In other words, the supremum norm of the function difference, denoted as $\left\Vert f - p \right\Vert_{\infty}$, is less than $\epsilon$~\cite{de1959stone}. Given a continuous real-valued function $f(x)$ and error bound $\epsilon>0$, we can apply Lagrange polynomial interpolation method based on Stone-Weierstrass Theorem to find an upper polynomial $p_{\mathscr{U}}(x) \geq f(x)$ and a lower polynomial $p_{\mathscr{L}}(x) \leq f(x)$ with respect to $x \in [m,M]$ such that:
\begin{eqnarray}\label{eq:lower and upper polynomials}
0 &\leq& p_{\mathscr{U}}(x) - f(x)~~\leq~~\epsilon, \nonumber \\
0 &\leq& f(x)-p_{\mathscr{L}}(x)~~\leq~~\epsilon,
\end{eqnarray}
In this paper, we also require (assumption) that $f(\bm{A})$ is a self-adjoint operator if $\bm{A}$ is a self-adjoint operator. 

Given two Hilbert spaces $\mathfrak{H}$ and $\mathfrak{K}$, $\mathbb{B}(\mathfrak{H})$ and $\mathbb{B}(\mathfrak{K})$ represent semi-algebra of all bounded linear operators on Hilbert spaces $\mathfrak{H}$ and $\mathfrak{K}$, respectively. Recall that $\Phi:\mathbb{B}(\mathfrak{H}) \rightarrow \mathbb{B}(\mathfrak{K})$ is a normalized positive linear map at the original settings of Choi-Davis-Jensen's Inequality. A normalized positive linear map is defined by the following Definition~\ref{def: normalized positive linear map}~\cite{pecaric2005mond,fujii2012recent}.
\begin{definition}\label{def: normalized positive linear map}
A map $\Phi:\mathbb{B}(\mathfrak{H}) \rightarrow \mathbb{B}(\mathfrak{K})$ is a normalized positive linear map if we have the three conditions satisfied jointly:
\begin{itemize}
\item The map $\Phi$ is a \emph{linear} map, i,e,  $\Phi(a\bm{X}+b\bm{Y})=a\Phi(\bm{X})+b\Phi(\bm{Y})$ for any $a,b \in \mathbb{C}$ and any $\bm{X}, \bm{Y} \in \mathbb{B}(\mathfrak{H})$.
\item The linear map $\Phi$ is a \emph{positive} map if the operator order (Loewner order) is preserved, i,e,  $\bm{X}\geq\bm{Y}$ implies $\Phi(\bm{X})\geq\Phi(\bm{Y})$.
\item The linear map $\Phi$ is a \emph{normalized} map if the identity operator is preserved, i,e,  $\Phi(\bm{I}_{\mathfrak{H}})=\bm{I}_{\mathfrak{K}}$, where $\bm{I}_{\mathfrak{H}}$ and $\bm{I}_{\mathfrak{K}}$ are identity operators of the Hilbert spaces $\mathfrak{H}$ and $\mathfrak{K}$, respectively. 
\end{itemize}
\end{definition}

In this work, we will consider more general $\Phi$ by assuming that 
\begin{eqnarray}\label{eq: new phi def}
\Phi(\bm{X})&=&\bm{V}^{\ast}\left(\sum\limits_{i=0}^{I}a_{i}\bm{X}^{i}\right)\bm{V}\nonumber \\
&=&\bm{V}^{\ast}\left(\sum\limits_{i_{+}\in S_{I_{+}}}a_{i_{+}}\bm{X}^{i_{+}}+\sum\limits_{i_{-}\in S_{I_{-}}}a_{i_{-}}\bm{X}^{i_{-}}\right)\bm{V},
\end{eqnarray}
where $\bm{V}$ is isometry in $\mathfrak{H}$ such that $\bm{V}^{\ast}\bm{V}=\bm{I}_{\mathfrak{H}}$, $a_{i_{+}}$ represent those nonnegative coefficients in $a_{i}$, and $a_{i_{-}}$ represent those negative coefficients in $a_{i}$. The collection of $i_{+}$ forms the set $S_{I_{+}}$, and the collection of $i_{-}$ forms the set $S_{I_{-}}$. Note that we do not require linear map, positive map, and normalized map for $\Phi$ defined by Eq.~\eqref{eq: new phi def}. Under the assumption provided by Eq.~\eqref{eq: new phi def}, the conventional normalized positive linear map provided by Definition~\ref{def: normalized positive linear map} is a special case by setting the polynomial $\sum\limits_{i=0}^{I}a_{i}\bm{X}^{i}$ as the identity map, i.e., all $a_i=0$ except $a_1=1$. 

We require the following Lemma~\ref{lma: lower and upper bound for f(A)} to provide the lower and upper bounds for $f(\bm{A})$ if $\bm{A}$ is a self-adjoint operator. 
\begin{lemma}\label{lma: lower and upper bound for f(A)}
Given a self-adjoint operator $\bm{A}$ with $\Lambda(\bm{A})$, such that 
\begin{eqnarray}\label{eq1: lma: lower and upper bound for f(A)}
0 &\leq& p_{\mathscr{U}}(x) - f(x)~~\leq~~\epsilon, \nonumber \\
0 &\leq& f(x)-p_{\mathscr{L}}(x)~~\leq~~\epsilon,
\end{eqnarray}
for $x \in [\min(\Lambda(\bm{A})), \max(\Lambda(\bm{A}))]$ with polynomials $p_{\mathscr{L}}(x)$ and $p_{\mathscr{U}}(x)$ expressed by 
\begin{eqnarray}\label{eq: lower poly formats}
p_{\mathscr{L}}(x)=\sum\limits_{k=0}^{n_{\mathscr{L}}}\alpha_k x^k, ~~
p_{\mathscr{U}}(x)=\sum\limits_{j=0}^{n_{\mathscr{U}}}\beta_j x^j. 
\end{eqnarray}
Further assume that $p_{\mathscr{L}}(\bm{A})\geq \bm{0}$, we have
\begin{eqnarray}\label{eq2-1: lma: lower and upper bound for f(A)}
\mathscr{K}^{-1}(\min(\Lambda(p_{\mathscr{L}}(\bm{A}))),\max(\Lambda(p_{\mathscr{L}}(\bm{A}))),i_+)p^{i_+}_{\mathscr{L}}(\bm{A})\leq f^{i_+}(\bm{A}),
\end{eqnarray}
and
\begin{eqnarray}\label{eq2-2: lma: lower and upper bound for f(A)}
f^{i_+}(\bm{A})\leq\mathscr{K}(\min(\Lambda(p_{\mathscr{U}}(\bm{A}))),\max(\Lambda(p_{\mathscr{U}}(\bm{A}))),i_+)p^{i_+}_{\mathscr{U}}(\bm{A}),
\end{eqnarray}
where Kantorovich functions $\mathscr{K}^{-1}(\min(\Lambda(p_{\mathscr{L}}(\bm{A}))),\max(\Lambda(p_{\mathscr{L}}(\bm{A}))),i_+)$ and \\$\mathscr{K}(\min(\Lambda(p_{\mathscr{U}}(\bm{A}))),\max(\Lambda(p_{\mathscr{U}}(\bm{A}))),i_+)$ are defined by Eq.~\eqref{eq: Kantorovich function}. Moreover, we also have 
\begin{eqnarray}\label{eq3-1: lma: lower and upper bound for f(A)}
\mathscr{K}(\min(\Lambda(p_{\mathscr{U}}(\bm{A}))),\max(\Lambda(p_{\mathscr{U}}(\bm{A}))),i_-)p^{i_-}_{\mathscr{U}}(\bm{A})\geq f^{i_-}(\bm{A}),
\end{eqnarray}
and
\begin{eqnarray}\label{eq3-2: lma: lower and upper bound for f(A)}
\mathscr{K}^{-1}(\min(\Lambda(p_{\mathscr{L}}(\bm{A}))),\max(\Lambda(p_{\mathscr{L}}(\bm{A}))),i_-)p^{i_-}_{\mathscr{L}}(\bm{A})\leq f^{i_-}(\bm{A}),
\end{eqnarray}
where Kantorovich function $\mathscr{K}(\min(\Lambda(p_{\mathscr{U}}(\bm{A}))),\max(\Lambda(p_{\mathscr{U}}(\bm{A}))),i_-)p^{i_-}_{\mathscr{U}}(\bm{A})$ and \\$\mathscr{K}^{-1}(\min(\Lambda(p_{\mathscr{L}}(\bm{A}))),\max(\Lambda(p_{\mathscr{L}}(\bm{A}))),i_-)p^{i_-}_{\mathscr{L}}(\bm{A})$ are defined by Eq.~\eqref{eq: Kantorovich function}.
\end{lemma}
\textbf{Proof:}
From Eq.~\eqref{eq1: lma: lower and upper bound for f(A)} and spectrum mapping theorem of the self-adjoint operator $\bm{A}$, we have 
\begin{eqnarray}\label{eq4: lma: lower and upper bound for f(A)}
f(\bm{A})\leq p_{\mathscr{U}}(\bm{A}),
\end{eqnarray}
and 
\begin{eqnarray}\label{eq5: lma: lower and upper bound for f(A)}
p_{\mathscr{L}}(\bm{A})\leq f(\bm{A}).
\end{eqnarray}
From Theorem 8.3 in~\cite{pecaric2005mond}, we have Eq.~\eqref{eq2-1: lma: lower and upper bound for f(A)} from Eq.~\eqref{eq5: lma: lower and upper bound for f(A)}. Similarly, we have Eq.~\eqref{eq2-2: lma: lower and upper bound for f(A)} from Eq.~\eqref{eq4: lma: lower and upper bound for f(A)}. 

Again, from Theorem 8.3 in~\cite{pecaric2005mond}, we have Eq.~\eqref{eq3-1: lma: lower and upper bound for f(A)} from Eq.~\eqref{eq4: lma: lower and upper bound for f(A)}. Similarly, we have Eq.~\eqref{eq3-2: lma: lower and upper bound for f(A)} from Eq.~\eqref{eq5: lma: lower and upper bound for f(A)}. 
$\hfill \Box$

\begin{remark}
Note that in Eq.~\eqref{eq2-1: lma: lower and upper bound for f(A)} and Eq.~\eqref{eq2-2: lma: lower and upper bound for f(A)}, we use the Kantorovich function \\
$\mathscr{K}^{-1}(\min(\Lambda(p_{\mathscr{L}}(\bm{A}))),\max(\Lambda(p_{\mathscr{L}}(\bm{A}))),i_+)$ and $\mathscr{K}(\min(\Lambda(p_{\mathscr{U}}(\bm{A}))),\max(\Lambda(p_{\mathscr{U}}(\bm{A}))),i_+)$ instead of using $\mathscr{K}(\min(\Lambda(f(\bm{A}))),\max(\Lambda(f(\bm{A}))),i_+)$ because $\Lambda(f(\bm{A}))$ may not within real domain. Same reason applies to Eq.~\eqref{eq3-1: lma: lower and upper bound for f(A)} and Eq.~\eqref{eq3-2: lma: lower and upper bound for f(A)}.
\end{remark}

Our next Lemma~\ref{lma: phi(f(A)) bounds} is about the upper and the lower bounds for $\Phi(f(\bm{A}))$ based on Lemma~\ref{lma: lower and upper bound for f(A)}.
\begin{lemma}\label{lma: phi(f(A)) bounds}
Under the definition of $\Phi$ provided by Eq.~\eqref{eq: new phi def} and same conditions provided by  Lemma~\ref{lma: lower and upper bound for f(A)}, we have
\begin{eqnarray}\label{eq1:lma: phi(f(A)) bounds}
\Phi(f(\bm{A})) &\leq& \bm{V}^{\ast}\left\{\sum\limits_{i_{+}\in S_{I_{+}}}a_{i_{+}}\mathscr{K}(\min(\Lambda(p_{\mathscr{U}}(\bm{A}))),\max(\Lambda(p_{\mathscr{U}}(\bm{A}))),i_+)p^{i_+}_{\mathscr{U}}(\bm{A})\right. \nonumber \\
&  &\left.+\sum\limits_{i_{-}\in S_{I_{-}}}a_{i_{-}}\mathscr{K}^{-1}(\min(\Lambda(p_{\mathscr{L}}(\bm{A}))),\max(\Lambda(p_{\mathscr{L}}(\bm{A}))),i_-)p^{i_-}_{\mathscr{L}}(\bm{A})\right\}\bm{V}.
\end{eqnarray}
On the other hand, we also have
\begin{eqnarray}\label{eq2:lma: phi(f(A)) bounds}
\Phi(f(\bm{A})) &\geq& \bm{V}^{\ast}\left\{\sum\limits_{i_{+}\in S_{I_{+}}}a_{i_{+}}\mathscr{K}^{-1}(\min(\Lambda(p_{\mathscr{L}}(\bm{A}))),\max(\Lambda(p_{\mathscr{L}}(\bm{A}))),i_+)p^{i_+}_{\mathscr{L}}(\bm{A})\right. \nonumber \\
&  &\left.+\sum\limits_{i_{-}\in S_{I_{-}}}a_{i_{-}}\mathscr{K}(\min(\Lambda(p_{\mathscr{U}}(\bm{A}))),\max(\Lambda(p_{\mathscr{U}}(\bm{A}))),i_-)p^{i_-}_{\mathscr{U}}(\bm{A})\right\}\bm{V}.
\end{eqnarray}
\end{lemma}
\textbf{Proof:}
This Lemma is proved by applying Lemma~\ref{lma: lower and upper bound for f(A)} to the definition of $\Phi$ provided by Eq.~\eqref{eq: new phi def}.
$\hfill \Box$

The following Lemma~\ref{lma: f(phi(A)) bounds} is to provide the upper and lower bounds for $f(\Phi(\bm{A}))$.
\begin{lemma}\label{lma: f(phi(A)) bounds}
Given the mapping $\Phi$ defined by Eq.~\eqref{eq: new phi def} and a self-adjoint operator $\bm{A}$, we will have the spectrum for $\Phi(\bm{A})$ denoted by $\Lambda(\Phi(\bm{A}))$.~\footnote{Note that the spectrum $\Lambda(\Phi(\bm{A}))$ is composed by real numbers only from $\Phi$ defined by Eq.~\eqref{eq: new phi def}.} From Stone-Weierstrass Theorem, we have
\begin{eqnarray}\label{eq1: lma: f(phi(A)) bounds}
0 &\leq& \tilde{p}_{\mathscr{U}}(x) - f(x)~~\leq~~\epsilon, \nonumber \\
0 &\leq& f(x)-\tilde{p}_{\mathscr{L}}(x)~~\leq~~\epsilon,
\end{eqnarray}
for $x \in [\min(\Lambda(\Phi(\bm{A}))), \max(\Lambda(\Phi(\bm{A})))]$. The polynomials $\tilde{p}_{\mathscr{L}}(x)$ and $\tilde{p}_{\mathscr{U}}(x)$ can be expressed by 
\begin{eqnarray}\label{eq2: lma: f(phi(A)) bounds}
\tilde{p}_{\mathscr{L}}(x)=\sum\limits_{k=0}^{\tilde{n}_{\mathscr{L}}}\tilde{\alpha}_k x^k, ~~~~
\tilde{p}_{\mathscr{U}}(x)=\sum\limits_{j=0}^{\tilde{n}_{\mathscr{U}}}\tilde{\beta}_j x^j. 
\end{eqnarray}
Then, we have 
\begin{eqnarray}\label{eq3: lma: f(phi(A)) bounds}
f(\Phi(\bm{A}))&\leq&\sum\limits_{j=0}^{\tilde{n}_{\mathscr{U}}}\tilde{\beta}_j \Phi^j(\bm{A}),
\end{eqnarray}
and
\begin{eqnarray}\label{eq4: lma: f(phi(A)) bounds}
f(\Phi(\bm{A}))&\geq&\sum\limits_{k=0}^{\tilde{n}_{\mathscr{L}}}\tilde{\alpha}_k \Phi^k(\bm{A}).
\end{eqnarray}
\end{lemma}
\textbf{Proof:}
From $\Phi$ defined by Eq.~\eqref{eq: new phi def}, the operator $\Phi(\bm{A})$ is a self-adjoint operator as the integer power of any self-adjoint operator $\bm{A}$ is a self-adjoint opeartor again. Then, we have Eq.~\eqref{eq3: lma: f(phi(A)) bounds} from the spectrum mapping theorem of the self-adjoint operator $\Phi(\bm{A})$ and the first condition provided by Eq.~\eqref{eq1: lma: f(phi(A)) bounds}. Similarly,  we also have Eq.~\eqref{eq4: lma: f(phi(A)) bounds} from the spectrum mapping theorem of the self-adjoint operator $\Phi(\bm{A})$ and the second condition provided by Eq.~\eqref{eq1: lma: f(phi(A)) bounds}. 
$\hfill \Box$

We are ready to present the main theorem of this work about generalized operator Jensen’s inequality without requirements of the function $f$ convexity and the mapping $\Phi$ to be a normalized positive linear map. 
\begin{theorem}\label{thm: gen Jensen}
Given a self-adjoint operator $\bm{A}$ with $\Lambda(\bm{A})$, such that 
\begin{eqnarray}\label{eq1: thm: gen Jensen}
0 &\leq& p_{\mathscr{U}}(x) - f(x)~~\leq~~\epsilon, \nonumber \\
0 &\leq& f(x)-p_{\mathscr{L}}(x)~~\leq~~\epsilon,
\end{eqnarray}
for $x \in [\min(\Lambda(\bm{A})), \max(\Lambda(\bm{A}))]$ with polynomials $p_{\mathscr{L}}(x)$ and $p_{\mathscr{U}}(x)$ expressed by 
\begin{eqnarray}\label{eq2: thm: gen Jensen}
p_{\mathscr{L}}(x)=\sum\limits_{k=0}^{n_{\mathscr{L}}}\alpha_k x^k, ~~
p_{\mathscr{U}}(x)=\sum\limits_{j=0}^{n_{\mathscr{U}}}\beta_j x^j. 
\end{eqnarray}
We assume that $p_{\mathscr{L}}(\bm{A})\geq \bm{0}$ and the mapping $\Phi:\mathbb{B}(\mathfrak{H}) \rightarrow \mathbb{B}(\mathfrak{K})$ is defined by Eq.~\eqref{eq: new phi def}. From such $\Phi$, we have
\begin{eqnarray}\label{eq3: thm: gen Jensen}
0 &\leq& \tilde{p}_{\mathscr{U}}(x) - f(x)~~\leq~~\epsilon, \nonumber \\
0 &\leq& f(x)-\tilde{p}_{\mathscr{L}}(x)~~\leq~~\epsilon,
\end{eqnarray}
for $x \in [\min(\Lambda(\Phi(\bm{A}))), \max(\Lambda(\Phi(\bm{A})))]$. The polynomals $\tilde{p}_{\mathscr{L}}(x)$ and $\tilde{p}_{\mathscr{U}}(x)$ can be expressed by 
\begin{eqnarray}\label{eq4: thm: gen Jensen}
\tilde{p}_{\mathscr{L}}(x)=\sum\limits_{k=0}^{\tilde{n}_{\mathscr{L}}}\tilde{\alpha}_k x^k, ~~~~
\tilde{p}_{\mathscr{U}}(x)=\sum\limits_{j=0}^{\tilde{n}_{\mathscr{U}}}\tilde{\beta}_j x^j. 
\end{eqnarray}
We use the following notation simplifications:
\begin{eqnarray}\label{eq5: thm: gen Jensen}
\frac{d}{e}\times W(\bm{A})&\define&\bm{V}^{\ast}\left\{\sum\limits_{i_{+}\in S_{I_{+}}}a_{i_{+}}\mathscr{K}(\min(\Lambda(p_{\mathscr{U}}(\bm{A}))),\max(\Lambda(p_{\mathscr{U}}(\bm{A}))),i_+)p^{i_+}_{\mathscr{U}}(\bm{A})\right. \nonumber \\
&  &\left.+\sum\limits_{i_{-}\in S_{I_{-}}}a_{i_{-}}\mathscr{K}^{-1}(\min(\Lambda(p_{\mathscr{L}}(\bm{A}))),\max(\Lambda(p_{\mathscr{L}}(\bm{A}))),i_-)p^{i_-}_{\mathscr{L}}(\bm{A})\right\}\bm{V},
\end{eqnarray}
where $d,e$ are any positive real numbers and $W(\bm{A})$ is an operator-valued function with $\bm{A}$ as its argument. 
\begin{eqnarray}\label{eq6: thm: gen Jensen}
\frac{c}{e}\times X(\bm{A})&\define&\bm{V}^{\ast}\left\{\sum\limits_{i_{+}\in S_{I_{+}}}a_{i_{+}}\mathscr{K}^{-1}(\min(\Lambda(p_{\mathscr{L}}(\bm{A}))),\max(\Lambda(p_{\mathscr{L}}(\bm{A}))),i_+)p^{i_+}_{\mathscr{L}}(\bm{A})\right. \nonumber \\
&  &\left.+\sum\limits_{i_{-}\in S_{I_{-}}}a_{i_{-}}\mathscr{K}(\min(\Lambda(p_{\mathscr{U}}(\bm{A}))),\max(\Lambda(p_{\mathscr{U}}(\bm{A}))),i_-)p^{i_-}_{\mathscr{U}}(\bm{A})\right\}\bm{V},
\end{eqnarray}
where $c,e$ are any positive real numbers and $X(\bm{A})$ is an operator-valued function with $\bm{A}$ as its argument. 
\begin{eqnarray}\label{eq7: thm: gen Jensen}
c\times Y(\bm{A})&\define&\sum\limits_{j=0}^{\tilde{n}_{\mathscr{U}}}\tilde{\beta}_j \Phi^j(\bm{A}),
\end{eqnarray}
where $Y(\bm{A})$ is an operator-valued function with $\bm{A}$ as its argument. 
\begin{eqnarray}\label{eq8: thm: gen Jensen}
d\times Z(\bm{A})&\define&\sum\limits_{k=0}^{\tilde{n}_{\mathscr{L}}}\tilde{\alpha}_k \Phi^k(\bm{A}), 
\end{eqnarray}
where $Z(\bm{A})$ is an operator-valued function with $\bm{A}$ as its argument. Then, we have
\begin{eqnarray}\label{eq9: thm: gen Jensen}
e\Phi(f(\bm{A})) - d(W(\bm{A})-Z(\bm{A})) \leq f(\Phi(\bm{A})) \leq e\Phi(f(\bm{A}))-c(X(\bm{A})-Y(\bm{A})).
\end{eqnarray}
\end{theorem}
\textbf{Proof:}
From Lemma~\ref{lma: phi(f(A)) bounds} and Eq.~\eqref{eq5: thm: gen Jensen} with  Eq.~\eqref{eq6: thm: gen Jensen}, we have
\begin{eqnarray}\label{eq10: thm: gen Jensen}
c\left(\frac{X(\bm{A})+Y(\bm{A})}{2}+\frac{X(\bm{A})-Y(\bm{A})}{2}\right)\leq e\Phi(f(\bm{A}))\leq d\left(\frac{W(\bm{A})+Z(\bm{A})}{2}+\frac{W(\bm{A})-Z(\bm{A})}{2}\right). 
\end{eqnarray}
Moreover, from Lemma~\ref{lma: f(phi(A)) bounds} and Eq.~\eqref{eq7: thm: gen Jensen} with  Eq.~\eqref{eq8: thm: gen Jensen}, we have
\begin{eqnarray}\label{eq11: thm: gen Jensen}
d\left(\frac{W(\bm{A})+Z(\bm{A})}{2}-\frac{W(\bm{A})-Z(\bm{A})}{2}\right)\leq f(\Phi(\bm{A})) \leq c\left(\frac{X(\bm{A})+Y(\bm{A})}{2}-\frac{X(\bm{A})-Y(\bm{A})}{2}\right).
\end{eqnarray}
Finally, this theorem is proved by rearrange related terms in Eq.~\eqref{eq10: thm: gen Jensen} and Eq.~\eqref{eq11: thm: gen Jensen}.
$\hfill \Box$
 
From Theorem~\ref{thm: gen Jensen}, we can have the following corollary about eigenvalues majorization inequality immediatedly.
\begin{corollary}\label{cor: majorization of thm: gen Jensen}
Given $\bm{A}$ as an $n \times n$ Hermitan matrix, and let $\left(\lambda_{\mathscr{L},1}, \lambda_{\mathscr{L},2}, \cdots, \lambda_{\mathscr{L},n}\right)$, $\left(\lambda_{1}, \lambda_{2}, \cdots, \lambda_{n}\right)$, and $\left(\lambda_{\mathscr{R},1}, \lambda_{\mathscr{R},2}, \cdots, \lambda_{\mathscr{R},n}\right)$, be eigenvalues of $e\Phi(f(\bm{A})) - d(W(\bm{A})-Z(\bm{A}))$, $f(\Phi(\bm{A}))$ and $e\Phi(f(\bm{A}))-c(X(\bm{A})-Y(\bm{A}))$ respectively. Then, we have
\begin{eqnarray}\label{eq1:cor: majorization of thm: gen Jensen}
\left(\lambda_{\mathscr{L},1}, \lambda_{\mathscr{L},2}, \cdots, \lambda_{\mathscr{L},n}\right)\preceq_{wk}
\left(\lambda_{1}, \lambda_{2}, \cdots, \lambda_{n}\right)\preceq_{wk}\left(\lambda_{\mathscr{R},1}, \lambda_{\mathscr{R},2}, \cdots, \lambda_{\mathscr{R},n}\right),
\end{eqnarray}
where $\preceq_{wk}$ represents the weakly majorization relation between entries of two vectors. 
\end{corollary}
\textbf{Proof:}
Given two $n \times n$ Hermitan matrices with $\bm{A}\leq\bm{B}$, from Loewner Theorem (1934)~\cite{pecaric2005mond,fujii2012recent}, we have
\begin{eqnarray}\label{eq2:cor: majorization of thm: gen Jensen}
\lambda_{\bm{A},i} \leq \lambda_{\bm{B},i},
\end{eqnarray}
where $\lambda_{\bm{A},i}$ and $\lambda_{\bm{B},i}$ are eigenvalues of Hermitan matrices $\bm{A}$ and $\bm{B}$, respectively. Therefore, Eq.~\eqref{eq2:cor: majorization of thm: gen Jensen} can be expressed by $\left(\lambda_{\bm{A},1},\lambda_{\bm{A},2},\cdots,\lambda_{\bm{A},n}\right)\preceq_{wk}\left(\lambda_{\bm{B},1},\lambda_{\bm{B},2},\cdots,\lambda_{\bm{B},n}\right)$. This corollary is proved by applying
Loewner Theorem (1934) to Theorem~\ref{thm: gen Jensen}. 
$\hfill \Box$

We will provide the following examples about Theorem~\ref{thm: gen Jensen}. Example~\ref{exp: special case of CDJ paper} will show that the work about Choi-Davis-Jensen's inequality without convexity studied in~\cite{micic2017choi} is a special case of Theorem~\ref{thm: gen Jensen}.
\begin{example}\label{exp: special case of CDJ paper}
Let us repeat Theorem 2.2 from~\cite{micic2017choi} here again with their notations for easy comparison with our work. If $\Phi$ is a normalized positive linear map, and $f:[m,M] \rightarrow (0, \infty)$ is a continuous twice differentiable function such that $\alpha \leq f''$ on $[m,M]$, we have
\begin{eqnarray}\label{eq1:exp: special case of CDJ paper}
\frac{1}{K(m,M,f)}\Phi(f(\bm{A}))+\frac{\alpha}{2K(m,M,f)}\left[(M+m)\Phi(\bm{A})-Mm-\Phi(\bm{A}^2)\right]
\nonumber \\
\leq f(\Phi(\bm{A}))~~~~~~~~~~~~~~~~~~~~~~~~~~~~~~~~~~~~~~~~~~~~~~~~~~~~~~~~~~~~~~~~~~~~~~~~~~~~~~~~~~~~~~~~\nonumber \\
\leq K(m,M,f)\Phi(f(\bm{A})) - \frac{\alpha}{2}\left[(M+m)\Phi(\bm{A})-Mm-\Phi(\bm{A})^2\right],
\end{eqnarray}
where $K(m,M,f)$ is defined by
\begin{eqnarray}\label{eq2:exp: special case of CDJ paper}
K(m,M,f)&=&\max\limits_{x \in [m,M]}\left(\frac{1}{f(x)}\left(\frac{M-x}{M-m}f(m)+\frac{x-m}{M-m}f(M)\right)\right).
\end{eqnarray}
Note that $K(m,M,f)$ is greater than 1. As $K(m,M,f)>1$, we have 
\begin{eqnarray}
\frac{1}{K(m,M,f)}\Phi(f(\bm{A})) \leq K(m,M,f)\Phi(f(\bm{A})),
\end{eqnarray}
therefore, the following inequality implies Eq.~\eqref{eq1:exp: special case of CDJ paper}:
\begin{eqnarray}\label{eq3:exp: special case of CDJ paper}
K(m,M,f)\Phi(f(\bm{A}))+\frac{\alpha}{2K(m,M,f)}\left[(M+m)\Phi(\bm{A})-Mm\bm{I}-\Phi(\bm{A})^2\right]
\nonumber \\
\leq f(\Phi(\bm{A}))~~~~~~~~~~~~~~~~~~~~~~~~~~~~~~~~~~~~~~~~~~~~~~~~~~~~~~~~~~~~~~~~~~~~~~~~~~~~~~~~~~~~~~~~\nonumber \\
\leq K(m,M,f)\Phi(f(\bm{A})) - \frac{\alpha}{2}\left[(M+m)\Phi(\bm{A})-Mm\bm{I}-\Phi(\bm{A})^2\right].
\end{eqnarray}

If we set the following parameters for Theorem~\ref{thm: gen Jensen} with $\alpha\geq 0$:
\begin{eqnarray}\label{eq4:exp: special case of CDJ paper}
c&=&\frac{\alpha}{2},\nonumber\\
d&=&\frac{\alpha}{2K(m,M,f)},\nonumber\\
e&=&K(m,M,f),\nonumber\\
W(\bm{A})&=&\Phi(\bm{A}^2)\geq_1\Phi(\bm{A})^2,\nonumber \\
X(\bm{A})&=&(M+m)\Phi(\bm{A})-Mm\bm{I},\nonumber \\
Y(\bm{A})&=&\Phi^2(\bm{A}),\nonumber \\
Z(\bm{A})&=&(M+m)\Phi(\bm{A})-Mm\bm{I}, 
\end{eqnarray}
where we use $\Phi$ as a normalized positive linear map in $\geq_1$. We will obtain the result provided by
Eq.~\eqref{eq3:exp: special case of CDJ paper} with $\alpha\geq 0$ which implies Eq.~\eqref{eq1:exp: special case of CDJ paper}. By comparing with the lower bounds of $f(\Phi(\bm{A}))$, our method can provide tighter lower bound compared to Theorem 2.2 from~\cite{micic2017choi}.

On the other hand, if $\alpha < 0$ and $\Phi(\bm{A})=\bm{V}^{\ast}\bm{A}\bm{V}$, we can set the following parameters to obtain  the result provided by Eq.~\eqref{eq2:exp: special case of CDJ paper} which implies Eq.~\eqref{eq1:exp: special case of CDJ paper}:
\begin{eqnarray}\label{eq4:exp: special case of CDJ paper}
c&=&\frac{-\alpha}{2},\nonumber\\
d&=&\frac{-\alpha}{2K(m,M,f)},\nonumber\\
e&=&K(m,M,f),\nonumber\\
W(\bm{A})&=&(M+m)\Phi(\bm{A})-Mm\bm{I},\nonumber \\
X(\bm{A})&=&\Phi^2(\bm{A}) =_1 \Phi(\bm{A}^2),\nonumber \\
Y(\bm{A})&=&(M+m)\Phi(\bm{A})-Mm\bm{I},\nonumber \\
Z(\bm{A})&=&\Phi(\bm{A})^2,
\end{eqnarray}
where we use $\Phi(\bm{A})=\bm{V}^{\ast}\bm{A}\bm{V}$ in $=_1$.
\end{example}

In the follownig example, we will consider the situation that $f$ is bounded by quadratic functions and $\Phi$ is also being generated by quadratic function. 
\begin{example}\label{exp: any func bounded by quadratic}
Given a self-adjoint operator $\bm{A}$ with $\Lambda(\bm{A})$, such that 
\begin{eqnarray}\label{eq1: exp: any func}
0 &\leq& p_{\mathscr{U}}(x) - f(x)~~\leq~~\epsilon, \nonumber \\
0 &\leq& f(x)-p_{\mathscr{L}}(x)~~\leq~~\epsilon,
\end{eqnarray}
for $x \in [\min(\Lambda(\bm{A})), \max(\Lambda(\bm{A}))]$ with polynomials $p_{\mathscr{L}}(x)$ and $p_{\mathscr{U}}(x)$ expressed by 
\begin{eqnarray}\label{eq2: exp: any func}
p_{\mathscr{L}}(x)=\sum\limits_{k=0}^{2}\alpha_k x^k, ~~
p_{\mathscr{U}}(x)=\sum\limits_{j=0}^{2}\beta_j x^j. 
\end{eqnarray}
We assume that $p_{\mathscr{L}}(\bm{A})\geq \bm{0}$ and the mapping $\Phi:\mathbb{B}(\mathfrak{H}) \rightarrow \mathbb{B}(\mathfrak{K})$ is defined as:
\begin{eqnarray}\label{eq2-1: exp: any func}
\Phi(\bm{A})&=&\bm{V}^{\ast}\left(a_0+a_1\bm{A}+a_2\bm{A}^2\right)\bm{V},
\end{eqnarray}
where we assume that $a_0, a_2 >0$ and $a_1 < 0$.

From such $\Phi$, we have
\begin{eqnarray}\label{eq3: exp: any func}
0 &\leq& \tilde{p}_{\mathscr{U}}(x) - f(x)~~\leq~~\epsilon, \nonumber \\
0 &\leq& f(x)-\tilde{p}_{\mathscr{L}}(x)~~\leq~~\epsilon,
\end{eqnarray}
for $x \in [\min(\Lambda(\Phi(\bm{A}))), \max(\Lambda(\Phi(\bm{A})))]$ with . The polynomals $\tilde{p}_{\mathscr{L}}(x)$ and $\tilde{p}_{\mathscr{U}}(x)$ are assumed to be 
\begin{eqnarray}\label{eq4: exp: any func}
\tilde{p}_{\mathscr{L}}(x)=p_{\mathscr{L}}(x), ~~~~
\tilde{p}_{\mathscr{U}}(x)=p_{\mathscr{U}}(x)
\end{eqnarray}

If we select real numbers $c,d,e$ as 1, we have the following notation simplifications:
\begin{eqnarray}\label{eq5: exp: any func}
W(\bm{A})&\define&\bm{V}^{\ast}\left\{a_0\bm{I}+ a_2\mathscr{K}(\min(\Lambda(p_{\mathscr{U}}(\bm{A}))),\max(\Lambda(p_{\mathscr{U}}(\bm{A}))),2)p^{2}_{\mathscr{U}}(\bm{A})\right. \nonumber \\
&  &\left.+a_{1}\mathscr{K}^{-1}(\min(\Lambda(p_{\mathscr{L}}(\bm{A}))),\max(\Lambda(p_{\mathscr{L}}(\bm{A}))),1)p_{\mathscr{L}}(\bm{A})\right\}\bm{V},
\end{eqnarray}
and
\begin{eqnarray}\label{eq6: exp: any func}
X(\bm{A})&\define&\bm{V}^{\ast}\left\{a_0\bm{I}+a_{2}\mathscr{K}^{-1}(\min(\Lambda(p_{\mathscr{L}}(\bm{A}))),\max(\Lambda(p_{\mathscr{L}}(\bm{A}))),2)p^{2}_{\mathscr{L}}(\bm{A})\right. \nonumber \\
&  &\left.+a_{1}\mathscr{K}(\min(\Lambda(p_{\mathscr{U}}(\bm{A}))),\max(\Lambda(p_{\mathscr{U}}(\bm{A}))),1)p^{i_-}_{\mathscr{U}}(\bm{A})\right\}\bm{V},
\end{eqnarray}
and
\begin{eqnarray}\label{eq7: exp: any func}
Y(\bm{A})&\define&\beta_0\bm{I}+\beta_1\bm{V}^{\ast}\left(a_0\bm{I}+a_1\bm{A}+a_2\bm{A}^2\right)\bm{V}\nonumber\\
&&+\beta_2\bm{V}^{\ast}\left(a^2_0\bm{I}+2a_0 a_1 \bm{A}+a^2_1\bm{A}^2+2a_0 a_2 \bm{A}^2+ 2a_1 a_2 \bm{A}^3+ a^2_2\bm{A}^4\right)\bm{V},
\end{eqnarray}
and
\begin{eqnarray}\label{eq8: exp: any func}
Z(\bm{A})&\define&\alpha_0\bm{I}+\alpha_1\bm{V}^{\ast}\left(a_0+a_1\bm{A}+a_2\bm{A}^2\right)\bm{V}\nonumber \\
&&+\alpha_2\bm{V}^{\ast}\left(a^2_0\bm{I}+2a_0 a_1 \bm{A}+a^2_1\bm{A}^2+2a_0 a_2 \bm{A}^2+ 2a_1 a_2 \bm{A}^3+ a^2_2\bm{A}^4\right)\bm{V}.
\end{eqnarray}
Then, we have
\begin{eqnarray}\label{eq9: exp: any func}
\Phi(f(\bm{A})) - (W(\bm{A})-Z(\bm{A})) \leq f(\Phi(\bm{A})) \leq \Phi(f(\bm{A}))-(X(\bm{A})-Y(\bm{A})).
\end{eqnarray}
\end{example}

\section{Applications}\label{sec: Applications}

In this section, we will apply proposed generalized Choi-Davis-Jensen's operator inequalities to statistical physics. Tsallis entropy for single random variable was introduced by C.Tsallis as \\ $T_q(X)= - \sum\limits_{x} p(x)^q f_q(p(x))$, incorporating a single parameter $q$ to extend the concept of Shannon entropy. Here, the $q$-logarithm is defined as $f_q(x)=\frac{x^{1-q}-1}{1-q}$ for any nonnegative real numbers $q$ and $x$, while $p(x)$ represents the probability distribution of the given random variable $X$. It is evident that the Tsallis entropy $T_q(X)$ converges to the Shannon entropy $-\sum\limits_{x} p(x) \log p(x)$ as $q$ approaches $1$, given that the $q$-logarithm uniformly converges to the natural logarithm as $q$ approaches 1. Tsallis entropy plays a crucial role in nonextensive statistics, often referred to as Tsallis statistics~\cite{furuichi2004fundamental}.

\subsection{Tsallis Relative Entropy Without $\Phi$}\label{sec: Tsallis Relative Entropy without Phi}

For Tsallis relative entropy for self-adjoint operators $\bm{A}>\bm{0}$ and $\bm{B}$, denoted by $T_q(\bm{A}\parallel \bm{B})$, it is defined as~\cite{micic2017choi}:
\begin{eqnarray}\label{eq: Tsallis relative entropy for operator}
T_q(\bm{A}\parallel \bm{B})&\define& \frac{\bm{A}\#_q\bm{B}-\bm{A}}{q},
\end{eqnarray}
where $-1 \leq q \leq 1$ with $q \neq 0$, and $\bm{A}\#_q\bm{B}\define\bm{A}^{1/2}\left(\bm{A}^{-1/2}\bm{B}\bm{A}^{-1/2}\right)^q\bm{A}^{1/2}$. We have the following Lemma~\ref{lma: Tsallis entropy bound no Phi} about $T_q(\bm{A}\parallel \bm{B})$ bounds.   
\begin{lemma}\label{lma: Tsallis entropy bound no Phi}
Given two positive self-adjoint operators $\bm{A}$ and $\bm{B}$ satisfying $m\bm{A} \leq \bm{B} \leq M\bm{A}$, where $0 < m <M$ with $m \geq 2$ and $M \geq 5m$, we define following relations:
\begin{eqnarray}\label{eq1: lma: Tsallis entropy bound no Phi}
\Gamma(\bm{A},\bm{B},q)&\define&-\frac{(\bm{B}-m\bm{A})(1-M^q)+(M\bm{A}-\bm{B})(1-m^q)}{q(M-m)},\nonumber\\
\Psi(\bm{A},\bm{B})&\define& \bm{A}\#_2\bm{B}-(M+m)\bm{B}+Mm\bm{A},\nonumber \\
\Omega(\bm{A},q)&\define&\frac{m^q\bm{A}}{(M+m)^{Mm}}.
\end{eqnarray}
Further, if we have $0< q \leq 1$ and $\bm{A}^{-1/2}\bm{B}\bm{A}^{-1/2}$ assumed to be self-adjoint, then, 
\begin{eqnarray}\label{eq1: lma: Tsallis entropy bound no Phi}
T_q(\bm{A}\parallel \bm{B})&\geq& \Gamma(\bm{A},\bm{B},q)-\frac{(1-q)M^{q-2}}{2}\Psi(\bm{A},\bm{B})+\Omega(\bm{A},q),\nonumber \\
T_q(\bm{A}\parallel \bm{B})&\leq& \Gamma(\bm{A},\bm{B},q)-\frac{(1-q)m^{q-2}}{2}\Psi(\bm{A},\bm{B})-\Omega(\bm{A},q).
\end{eqnarray}
\end{lemma}
\textbf{Proof:}
Because we have the following inequalites with respect to the function $\frac{1-x^{q}}{q}$ given $0 < m <M$ with $m \geq 2$ and $M \geq 5m$:
\begin{eqnarray}\label{eq2: lma: Tsallis entropy bound no Phi}
\frac{1-x^{q}}{q}&\geq&\frac{(x-m)(1-M^q)+(M-x)(1-m^q)}{q(M-m)}+\frac{(1-q)m^{q-2}}{2}\left(x^2-(M+m)x+Mm\right)\nonumber \\
&&+\frac{m^q}{(M+m)^{Mm}},\nonumber \\
\frac{1-x^{q}}{q}&\leq&\frac{(x-m)(1-M^q)+(M-x)(1-m^q)}{q(M-m)}+\frac{(1-q)M^{q-2}}{2}\left(x^2-(M+m)x+Mm\right)\nonumber \\
&&-\frac{m^q}{(M+m)^{Mm}}.
\end{eqnarray}
This lemma is proved by setting $x$ with the positive self-adjoint operator $\bm{A}^{-1/2}\bm{B}\bm{A}^{-1/2}$ in Eq.~\eqref{eq2: lma: Tsallis entropy bound no Phi} followed by multiplying both sides by $\bm{A}^{1/2}$.
$\hfill \Box$

From Lemma~\ref{lma: Tsallis entropy bound no Phi}, we have the following lemma immediately about tighter bounds of relative operator entropy provided by~\cite{micic2017choi} by taking $q \rightarrow 0$. 
\begin{lemma}\label{lma: bounds of relative operator entropy}
The relative operator entropy with respect to operators $\bm{A}$ and $\bm{B}$, denoted by $S(\bm{A}\parallel \bm{B})$, is given by
\begin{eqnarray}\label{eq1: lma: bounds of relative operator entropy}
S(\bm{A}\parallel \bm{B})&\define&\bm{A}^{1/2}\log\left(\bm{A}^{-1/2}\bm{B}\bm{A}^{-1/2}\right)\bm{A}^{1/2}.
\end{eqnarray}
Under same conditions provided by Lemma~\ref{lma: Tsallis entropy bound no Phi}, we have
\begin{eqnarray}\label{eq2: lma: bounds of relative operator entropy}
S(\bm{A}\parallel \bm{B})&\geq&\frac{(\bm{B}-m\bm{A})\log M+(M\bm{A}-\bm{B})\log m}{(M-m)}-\frac{1}{2M^2}\Psi(\bm{A},\bm{B})+\Omega(\bm{A},0),\nonumber \\
S(\bm{A}\parallel \bm{B})&\leq&\frac{(\bm{B}-m\bm{A})\log M+(M\bm{A}-\bm{B})\log m}{(M-m)}-\frac{1}{2m^2}\Psi(\bm{A},\bm{B})-\Omega(\bm{A},0).
\end{eqnarray}
\end{lemma}
\textbf{Proof:}
Because we have
\begin{eqnarray}\label{eq3: lma: bounds of relative operator entropy}
\lim\limits_{q \rightarrow 0} \frac{1-x^{q}}{q} = - \log x,
\end{eqnarray}  
the term $\Gamma(\bm{A},\bm{B},q)$ becomes
\begin{eqnarray}\label{eq3: lma: bounds of relative operator entropy}
\lim\limits_{q \rightarrow 0} \Gamma(\bm{A},\bm{B},q)
&=&\frac{(\bm{B}-m\bm{A})\log M+(M\bm{A}-\bm{B})\log m}{(M-m)}.
\end{eqnarray}
This lemma is proved by applying Eq.~\eqref{eq3: lma: bounds of relative operator entropy} to Lemma~\ref{lma: Tsallis entropy bound no Phi}.
$\hfill \Box$

\begin{remark}
Lemma~\ref{lma: bounds of relative operator entropy} provides tighter bounds of relative operator entropy given by  by~\cite{micic2017choi} as $\Omega(\bm{A},0)$ is a positive operator. 
\end{remark}

\subsection{Tsallis Relative Entropy With $\Phi$}\label{sec: Tsallis Relative Entropy with Phi}

In section~\ref{sec: Tsallis Relative Entropy without Phi}, operators $\bm{A}$ and $\bm{B}$ in Tsallis relative entropy are applied directly into Eq.~\eqref{eq: Tsallis relative entropy for operator}. In this section, we will consider applying $\Phi$ to operators $\bm{A}$, $\bm{B}$ and $T_q(\bm{A}\parallel \bm{B})$. The mapping $\Phi$ is assumed to have the following format:
\begin{eqnarray}\label{eq: new phi def in app}
\Phi(\bm{X})&=&\bm{V}^{\ast}\left(a_{0}\bm{I}+a_1\bm{X}+a_2\bm{X}^2\right)\bm{V},
\end{eqnarray}
where coefficients $a_0,a_1$ and $a_2$ are positive real numbers. We have the following Lemma~\ref{lma: Tsallis entropy bound Phi} about Tsallis relative entropy with respect to $\Phi(\bm{A})$ and $\Phi(\bm{B})$. 

\begin{lemma}\label{lma: Tsallis entropy bound Phi}
Given two positive self-adjoint operators $\bm{A}$ and $\bm{B}$ satisfying $m\Phi(\bm{A}) \leq \Phi(\bm{B}) \leq M\Phi(\bm{A})$, where $\Phi$ is given by Eq.~\eqref{eq: new phi def in app} and constants $M,m$ satisfy $0 < m <M$ with $m \geq 2$ and $M \geq 5m$, we define following relations:
\begin{eqnarray}\label{eq1: lma: Tsallis entropy bound Phi}
\Gamma(\Phi(\bm{A}),\Phi(\bm{B}),q)&\define&-\frac{(\Phi(\bm{B})-m\Phi(\bm{A}))(1-M^q)+(M\Phi(\bm{A})-\Phi(\bm{B}))(1-m^q)}{q(M-m)},\nonumber\\
\Psi(\Phi(\bm{A}),\Phi(\bm{B}))&\define& \Phi(\bm{A})\#_2\Phi(\bm{B})-(M+m)\Phi(\bm{B})+Mm\Phi(\bm{A}),\nonumber \\
\Omega(\Phi(\bm{A}),q)&\define&\frac{m^q\Phi(\bm{A})}{(M+m)^{Mm}}.
\end{eqnarray}
Further, if we have $0< q \leq 1$, then, 
\begin{eqnarray}\label{eq1: lma: Tsallis entropy bound Phi}
T_q(\Phi(\bm{A})\parallel\Phi(\bm{B}))&\geq& \Gamma(\Phi(\bm{A}),\Phi(\bm{B}),q)-\frac{(1-q)M^{q-2}}{2}\Psi(\Phi(\bm{A}),\Phi(\bm{B}))+\Omega(\Phi(\bm{A}),q),\nonumber \\
T_q(\Phi(\bm{A})\parallel\Phi(\bm{B}))&\leq& \Gamma(\Phi(\bm{A}),\Phi(\bm{B}),q)-\frac{(1-q)m^{q-2}}{2}\Psi(\Phi(\bm{A}),\Phi(\bm{B}))-\Omega(\Phi(\bm{A}),q).
\end{eqnarray}
\end{lemma}
\textbf{Proof:}
Since we have the following inequalites with respect to the function $\frac{1-x^{q}}{q}$ given $0 < m <M$ with $m \geq 2$ and $M \geq 5m$:
\begin{eqnarray}\label{eq2: lma: Tsallis entropy bound Phi}
\frac{1-x^{q}}{q}&\geq&\frac{(x-m)(1-M^q)+(M-x)(1-m^q)}{q(M-m)}+\frac{(1-q)m^{q-2}}{2}\left(x^2-(M+m)x+Mm\right)\nonumber \\
&&+\frac{m^q}{(M+m)^{Mm}},\nonumber \\
\frac{1-x^{q}}{q}&\leq&\frac{(x-m)(1-M^q)+(M-x)(1-m^q)}{q(M-m)}+\frac{(1-q)M^{q-2}}{2}\left(x^2-(M+m)x+Mm\right)\nonumber \\
&&-\frac{m^q}{(M+m)^{Mm}}.
\end{eqnarray}
This lemma is proved by setting $x$ with the positive self-adjoint operator $\Phi(\bm{A})^{-1/2}\Phi(\bm{B})\Phi(\bm{A})^{-1/2}$ in Eq.~\eqref{eq2: lma: Tsallis entropy bound Phi} followed by multiplying both sides by $\bm{A}^{1/2}$.
$\hfill \Box$

Lemma~\ref{lma: Phi of Tsallis entropy bound} below is established to provide bounds for $\Phi(T_q(\bm{A}\parallel\bm{B}))$.
\begin{lemma}\label{lma: Phi of Tsallis entropy bound}
Given $\Phi$ defined by Eq.~\eqref{eq: new phi def in app}, and two positive self-adjoint operators $\bm{A}$ and $\bm{B}$ satisfying $m\bm{A} \leq \bm{B} \leq M\bm{A}$, where $0 < m <M$ with $m \geq 2$ and $M \geq 5m$, we set up the follownig two polynomials $p_{\mathscr{L}}(x)$ and $p_{\mathscr{U}}(x)$:
\begin{eqnarray}\label{eq1: lma: Phi of Tsallis entropy bound}
p_{\mathscr{L}}(x)&=&-\frac{(x-m)(1-M^q)+(M-x)(1-m^q)}{q(M-m)}-\frac{(1-q)M^{q-2}}{2}\left(x^2-(M+m)x+Mm\right)\nonumber \\
&&+\frac{m^q}{(M+m)^{Mm}},\nonumber \\
p_{\mathscr{U}}(x)&=&-\frac{(x-m)(1-M^q)+(M-x)(1-m^q)}{q(M-m)}-\frac{(1-q)m^{q-2}}{2}\left(x^2-(M+m)x+Mm\right)\nonumber \\
&&-\frac{m^q}{(M+m)^{Mm}},
\end{eqnarray}
where $0 < q \leq 1$. Moreover, we also assume that $\bm{A}^{1/2}$ and $(\bm{A}^{-1/2}\bm{B}\bm{A}^{-1/2})^q$ are commute, then, we have
\begin{eqnarray}\label{eq2-1: lma: Phi of Tsallis entropy bound}
\Phi(T_q(\bm{A}\parallel\bm{B}))&\geq&\bm{V}^{\ast}\left\{a_0\bm{I}+ a_1\mathscr{K}^{-1}(\min(\Lambda(p_{\mathscr{L}}(\bm{A}^{-1/2}\bm{B}\bm{A}^{-1/2}))),\max(\Lambda(p_{\mathscr{L}}(\bm{A}^{-1/2}\bm{B}\bm{A}^{-1/2}))),1)\right.\nonumber \\
&  &\left.\cdot\bm{A}^{1/2}p_{\mathscr{L}}(\bm{A}^{-1/2}\bm{B}\bm{A}^{-1/2})\bm{A}^{1/2}\right. \nonumber \\
&  &\left.+a_{2}\mathscr{K}^{-1}(\min(\Lambda(p_{\mathscr{L}}(\bm{A}^{-1/2}\bm{B}\bm{A}^{-1/2}))),\max(\Lambda(p_{\mathscr{L}}(\bm{A}^{-1/2}\bm{B}\bm{A}^{-1/2}))),2)\right.\nonumber \\
&  &\left.\cdot\bm{A}p^2_{\mathscr{L}}(\bm{A}^{-1/2}\bm{B}\bm{A}^{-1/2})\bm{A}\right\}\bm{V}, \end{eqnarray}
\begin{eqnarray}\label{eq2-2: lma: Phi of Tsallis entropy bound}
\Phi(T_q(\bm{A}\parallel\bm{B}))&\leq&\bm{V}^{\ast}\left\{a_0\bm{I}+ a_1\mathscr{K}(\min(\Lambda(p_{\mathscr{U}}(\bm{A}^{-1/2}\bm{B}\bm{A}^{-1/2}))),\max(\Lambda(p_{\mathscr{U}}(\bm{A}^{-1/2}\bm{B}\bm{A}^{-1/2}))),1)\right.\nonumber \\
&  &\left.\cdot\bm{A}^{1/2}p_{\mathscr{U}}(\bm{A}^{-1/2}\bm{B}\bm{A}^{-1/2})\bm{A}^{1/2}\right. \nonumber \\
&  &\left.+a_{2}\mathscr{K}(\min(\Lambda(p_{\mathscr{U}}(\bm{A}^{-1/2}\bm{B}\bm{A}^{-1/2}))),\max(\Lambda(p_{\mathscr{U}}(\bm{A}^{-1/2}\bm{B}\bm{A}^{-1/2}))),2)\right.\nonumber \\
&  &\left.\cdot\bm{A}p^2_{\mathscr{U}}(\bm{A}^{-1/2}\bm{B}\bm{A}^{-1/2})\bm{A}\right\}\bm{V}. 
\end{eqnarray}
\end{lemma}
\textbf{Proof:}
Because we have
\begin{eqnarray}\label{eq3: lma: Phi of Tsallis entropy bound}
\frac{x^{q}-1}{q}&\geq&-\frac{(x-m)(1-M^q)+(M-x)(1-m^q)}{q(M-m)}-\frac{(1-q)M^{q-2}}{2}\left(x^2-(M+m)x+Mm\right)\nonumber \\
&&+\frac{m^q}{(M+m)^{Mm}}=_1 p_{\mathscr{L}}(x),
\end{eqnarray}
where $=_1$ comes from Eq.~\eqref{eq1: lma: Phi of Tsallis entropy bound}, by applying $\bm{A}^{-1/2}\bm{B}\bm{A}^{-1/2}$ to $x$ in Eq.~\eqref{eq3: lma: Phi of Tsallis entropy bound}, we also have
\begin{eqnarray}\label{eq6: lma: Phi of Tsallis entropy bound}
\Phi\left(T_q(\bm{A}\parallel\bm{B})\right)&=&\Phi\left(\bm{A}^{1/2}\left(\frac{(\bm{A}^{-1/2}\bm{B}\bm{A}^{-1/2})^q-\bm{I}}{q}\right)\bm{A}^{1/2}\right) \nonumber \\
&=_1& \bm{V}^{\ast}\left\{a_0 \bm{I} + a_1 \bm{A}^{1/2}\left(\frac{(\bm{A}^{-1/2}\bm{B}\bm{A}^{-1/2})^q-\bm{I}}{q}\right)\bm{A}^{1/2} \right.\nonumber \\
& &+\left.
a_2 \bm{A}\left(\frac{(\bm{A}^{-1/2}\bm{B}\bm{A}^{-1/2})^q-\bm{I}}{q}\right)^2\bm{A}
\right\}\bm{V}.\nonumber \\
&\geq_2&\bm{V}^{\ast}\left\{a_0\bm{I}+ a_1\mathscr{K}^{-1}(\min(\Lambda(p_{\mathscr{L}}(\bm{A}^{-1/2}\bm{B}\bm{A}^{-1/2}))),\max(\Lambda(p_{\mathscr{L}}(\bm{A}^{-1/2}\bm{B}\bm{A}^{-1/2}))),1)\right.\nonumber \\
&  &\left.\cdot \bm{A}^{1/2}p_{\mathscr{L}}(\bm{A}^{-1/2}\bm{B}\bm{A}^{-1/2})\bm{A}^{1/2}\right. \nonumber \\
&  &\left.+a_{2}\mathscr{K}^{-1}(\min(\Lambda(p_{\mathscr{L}}(\bm{A}^{-1/2}\bm{B}\bm{A}^{-1/2}))),\max(\Lambda(p_{\mathscr{L}}(\bm{A}^{-1/2}\bm{B}\bm{A}^{-1/2}))),2)\right.\nonumber \\
&  &\left.\cdot \bm{A}p^2_{\mathscr{L}}(\bm{A}^{-1/2}\bm{B}\bm{A}^{-1/2})\bm{A}\right\}\bm{V},
\end{eqnarray}
where we use that $\bm{A}^{1/2}$ and $(\bm{A}^{-1/2}\bm{B}\bm{A}^{-1/2})^q$ are commutative in $=_1$, and Lemma~\ref{lma: phi(f(A)) bounds} is used to obtain $\geq_2$. This establish Eq.~\eqref{eq2-1: lma: Phi of Tsallis entropy bound}.

Similarly, we can apply the following inequality to get Eq.~\eqref{eq2-2: lma: Phi of Tsallis entropy bound} by the same argument to derive Eq.~\eqref{eq2-1: lma: Phi of Tsallis entropy bound}:
\begin{eqnarray}
\frac{x^{q}-1}{q}&\leq&-\frac{(x-m)(1-M^q)+(M-x)(1-m^q)}{q(M-m)}-\frac{(1-q)m^{q-2}}{2}\left(x^2-(M+m)x+Mm\right)\nonumber \\
&&-\frac{m^q}{(M+m)^{Mm}}.
\end{eqnarray}
$\hfill \Box$

We are ready to present generalized Choi-Davis-Jensen's operator inequalities for Tsallis relative entropy. 
\begin{theorem}\label{thm: Tsallis entropy bound}
As conditions provided by Lemma~\ref{lma: Tsallis entropy bound Phi} and Lemma~\ref{lma: Phi of Tsallis entropy bound}, we have the following inequalities
\begin{eqnarray}\label{eq1: thm: Tsallis entropy bound}
\Phi(T_q(\bm{A}\parallel\bm{B})) - (W(\bm{A},\bm{B})-Z(\bm{A},\bm{B}))&\leq&T_q(\Phi(\bm{A})\parallel\Phi(\bm{B})) \nonumber \\
&\leq&\Phi(T_q(\bm{A}\parallel\bm{B}))-(X(\bm{A},\bm{B})-Y(\bm{A},\bm{B})).
\end{eqnarray}
where we have
\begin{eqnarray}
W(\bm{A},\bm{B})&=&\bm{V}^{\ast}\left\{a_0\bm{I}+ a_1\mathscr{K}(\min(\Lambda(p_{\mathscr{U}}(\bm{A}^{-1/2}\bm{B}\bm{A}^{-1/2}))),\max(\Lambda(p_{\mathscr{U}}(\bm{A}^{-1/2}\bm{B}\bm{A}^{-1/2}))),1)\right.\nonumber \\
&  &\left.\cdot\bm{A}^{1/2}p_{\mathscr{L}}(\bm{A}^{-1/2}\bm{B}\bm{A}^{-1/2})\bm{A}^{1/2}\right. \nonumber \\
&  &\left.+a_{2}\mathscr{K}(\min(\Lambda(p_{\mathscr{U}}(\bm{A}^{-1/2}\bm{B}\bm{A}^{-1/2}))),\max(\Lambda(p_{\mathscr{U}}(\bm{A}^{-1/2}\bm{B}\bm{A}^{-1/2}))),2)\right.\nonumber \\
&  &\left.\cdot\bm{A}p^2_{\mathscr{L}}(\bm{A}^{-1/2}\bm{B}\bm{A}^{-1/2})\bm{A}\right\}\bm{V},
\end{eqnarray}
\begin{eqnarray}
X(\bm{A},\bm{B})&=&\bm{V}^{\ast}\left\{a_0\bm{I}+ a_1\mathscr{K}^{-1}(\min(\Lambda(p_{\mathscr{L}}(\bm{A}^{-1/2}\bm{B}\bm{A}^{-1/2}))),\max(\Lambda(p_{\mathscr{L}}(\bm{A}^{-1/2}\bm{B}\bm{A}^{-1/2}))),1)\right.\nonumber \\
&  &\left.\cdot\bm{A}^{1/2}p_{\mathscr{L}}(\bm{A}^{-1/2}\bm{B}\bm{A}^{-1/2})\bm{A}^{1/2}\right. \nonumber \\
&  &\left.+a_{2}\mathscr{K}^{-1}(\min(\Lambda(p_{\mathscr{L}}(\bm{A}^{-1/2}\bm{B}\bm{A}^{-1/2}))),\max(\Lambda(p_{\mathscr{L}}(\bm{A}^{-1/2}\bm{B}\bm{A}^{-1/2}))),2)\right.\nonumber \\
&  &\left.\cdot\bm{A}p^2_{\mathscr{L}}(\bm{A}^{-1/2}\bm{B}\bm{A}^{-1/2})\bm{A}\right\}\bm{V},
\end{eqnarray}
\begin{eqnarray}
Y(\bm{A},\bm{B})&=&\Gamma(\Phi(\bm{A}),\Phi(\bm{B}),q)-\frac{(1-q)m^{q-2}}{2}\Psi(\Phi(\bm{A}),\Phi(\bm{B}))-\Omega(\Phi(\bm{A}),q),
\end{eqnarray}
\begin{eqnarray}
Z(\bm{A},\bm{B})&=& \Gamma(\Phi(\bm{A}),\Phi(\bm{B}),q)-\frac{(1-q)M^{q-2}}{2}\Psi(\Phi(\bm{A}),\Phi(\bm{B}))+\Omega(\Phi(\bm{A}),q).
\end{eqnarray}
\end{theorem}  
\textbf{Proof:}
By setting $c=d=e=1$ and replacing $W(\bm{A}),X(\bm{A}),Y(\bm{A})$ and $W(\bm{A})$ with $W(\bm{A},\bm{B}),X(\bm{A},\bm{B})$, $Y(\bm{A},\bm{B})$ and $Z(\bm{A},\bm{B})$, respectively, obtained from Lemma~\ref{lma: Tsallis entropy bound Phi} and Lemma~\ref{lma: Phi of Tsallis entropy bound} in Theorem~\ref{thm: gen Jensen}.
$\hfill \Box$

\begin{remark}
All results discussed in this work can also be applied to tensors if tensors under Einstein products are treated as operators~\cite{chang2023tail,chang2022random}.
\end{remark}

\bibliographystyle{IEEETran}
\bibliography{Gen_Jensen_Bib}

\begin{thebibliography}{10}
\providecommand{\url}[1]{#1}
\csname url@samestyle\endcsname
\providecommand{\newblock}{\relax}
\providecommand{\bibinfo}[2]{#2}
\providecommand{\BIBentrySTDinterwordspacing}{\spaceskip=0pt\relax}
\providecommand{\BIBentryALTinterwordstretchfactor}{4}
\providecommand{\BIBentryALTinterwordspacing}{\spaceskip=\fontdimen2\font plus
\BIBentryALTinterwordstretchfactor\fontdimen3\font minus
  \fontdimen4\font\relax}
\providecommand{\BIBforeignlanguage}[2]{{%
\expandafter\ifx\csname l@#1\endcsname\relax
\typeout{** WARNING: IEEEtran.bst: No hyphenation pattern has been}%
\typeout{** loaded for the language `#1'. Using the pattern for}%
\typeout{** the default language instead.}%
\else
\language=\csname l@#1\endcsname
\fi
#2}}
\providecommand{\BIBdecl}{\relax}
\BIBdecl

\bibitem{choi1974schwarz}
M.-D. Choi, ``A schwarz inequality for positive linear maps on
  $\bm{C}^{\ast}$-algebras,'' \emph{Illinois Journal of Mathematics}, vol.~18,
  no.~4, pp. 565--574, 1974.

\bibitem{davis1957schwarz}
C.~Davis, ``A schwarz inequality for convex operator functions,''
  \emph{Proceedings of the American Mathematical Society}, vol.~8, no.~1, pp.
  42--44, 1957.

\bibitem{becker1989functional}
L.~C. Becker, T.~Burton, and S.~Zhang, ``Functional differential equations and
  {J}ensen's inequality,'' \emph{Journal of Mathematical Analysis and
  Applications}, vol. 138, no.~1, pp. 137--156, 1989.

\bibitem{briet2009properties}
J.~Bri{\"e}t and P.~Harremo{\"e}s, ``Properties of classical and quantum
  {J}ensen-{S}hannon divergence,'' \emph{Physical Review A}, vol.~79, no.~5, p.
  052311, 2009.

\bibitem{briat2011convergence}
C.~Briat, ``Convergence and equivalence results for the {J}ensen's
  inequality—application to time-delay and sampled-data systems,'' \emph{IEEE
  Transactions on Automatic Control}, vol.~56, no.~7, pp. 1660--1665, 2011.

\bibitem{majtey2005jensen}
A.~P. Majtey, P.~W. Lamberti, and D.~P. Prato, ``{J}ensen-{S}hannon divergence
  as a measure of distinguishability between mixed quantum states,''
  \emph{Physical Review A}, vol.~72, no.~5, p. 052310, 2005.

\bibitem{futami2021loss}
F.~Futami, T.~Iwata, I.~Sato, M.~Sugiyama \emph{et~al.}, ``Loss function based
  second-order {J}ensen inequality and its application to particle variational
  inference,'' \emph{Advances in Neural Information Processing Systems},
  vol.~34, pp. 6803--6815, 2021.

\bibitem{chang2023algebraic}
S.-Y. Chang, ``Algebraic connectivity characterization of ensemble random
  hypergraphs,'' \emph{arXiv preprint arXiv:2310.08700}, 2023.

\bibitem{chang2023tail}
------, ``Tail bounds for multivariate random tensor means,'' \emph{arXiv
  preprint arXiv:2308.06478}, 2023.

\bibitem{chang2023random}
------, ``Random tensor inequalities and tail bounds for bivariate random
  tensor means, part i,'' \emph{arXiv preprint arXiv:2305.03301}, 2023.

\bibitem{chang2023tailRP}
------, ``Tail bounds for tensor-valued random process,'' \emph{arXiv preprint
  arXiv:2302.00602}, 2023.

\bibitem{chang2022general}
S.~Y. Chang and Y.~Wei, ``General tail bounds for random tensors summation:
  majorization approach,'' \emph{Journal of Computational and Applied
  Mathematics}, vol. 416, p. 114533, 2022.

\bibitem{chang2022random}
S.-Y. Chang, ``Random multiple operator integrals,'' \emph{arXiv preprint
  arXiv:2210.09392}, 2022.

\bibitem{chang2022convenient}
S.~Y. Chang and W.-W. Lin, ``Convenient tail bounds for sums of random
  tensors,'' \emph{Taiwanese Journal of Mathematics}, vol.~26, no.~3, pp.
  571--606, 2022.

\bibitem{micic2017choi}
J.~Mi{\'c}i{\'c}, H.~R. Moradi, and S.~Furuichi, ``{C}hoi-{D}avis-{J}ensen's
  inequality without convexity,'' \emph{arXiv preprint arXiv:1705.09784}, 2017.

\bibitem{hashemi2021reversing}
S.~S. Hashemi~Karouei, M.~S. Asgari, and M.~Shah~Hosseini, ``On reversing
  operator choi--davis--jensen inequality,'' \emph{Iranian Journal of Science
  and Technology, Transactions A: Science}, vol.~45, no.~4, pp. 1405--1410,
  2021.

\bibitem{de1959stone}
L.~De~Branges, ``The {S}tone-{W}eierstrass theorem,'' \emph{Proceedings of the
  American Mathematical Society}, vol.~10, no.~5, pp. 822--824, 1959.

\bibitem{pecaric2005mond}
J.~Pecaric, T.~Furuta, J.~M. Hot, and Y.~Seo, \emph{Mond-Pecaric method in
  operator inequalities}.\hskip 1em plus 0.5em minus 0.4em\relax Element
  Zagreb, Croatia, 2005.

\bibitem{fujii2012recent}
M.~Fujii and J.~M. Hot, ``Recent developments of mond-pecaric method in
  operator inequalities,'' \emph{Monographs in Inequalities}, vol.~4, 2012.

\bibitem{furuichi2004fundamental}
S.~Furuichi, K.~Yanagi, and K.~Kuriyama, ``Fundamental properties of {T}sallis
  relative entropy,'' \emph{Journal of Mathematical Physics}, vol.~45, no.~12,
  pp. 4868--4877, 2004.

\end{thebibliography}

\end{document}